\documentclass[12pt]{article}

\usepackage{amsmath}
\usepackage{natbib}
\setcitestyle{numbers}
\setcitestyle{super}
\setcitestyle{square}
\setcitestyle{sort}
\setcitestyle{comma}

\usepackage{times}
\usepackage{graphicx}
\graphicspath{ {.} }

\usepackage{placeins}

\let\Oldsection\section
\renewcommand{\section}{\FloatBarrier\Oldsection}

\let\Oldsubsection\subsection
\renewcommand{\subsection}{\FloatBarrier\Oldsubsection}

\let\Oldsubsubsection\subsubsection
\renewcommand{\subsubsection}{\FloatBarrier\Oldsubsubsection}

\newenvironment{sciabstract}{%
\begin{quote} \bf}
{\end{quote}}

\title{Ellipse Combining with Unknown Cross Ellipse Correlations}

\author
{Adam Hall}

\date{}

\begin{document} 


\baselineskip24pt


\maketitle


\begin{sciabstract}
  We discuss the combining of measurements where single measurement covariances are given but the joint measurement covariance is unknown. For this paper we assume the mapping of a single measurement to the solution space is the identity matrix. We examine the solution when it is assumed all measurements are uncorrelated. We then present a way to parameter joint measurement covariance based on pairwise correlation coefficients. Finally, we discuss how to use this parameterization to combine the measurements.
\end{sciabstract}


\section*{Introduction}

In this paper, we are going to examine the linear combining or fusion of location estimates. Each estimate is a vector of size $k\times 1$ and for each estimate we are also provided a $k\times k$ covariance matrix. We will assume that the overall joint system covariance which describes the correlation between location estimates is not provided. Specifically, we consider the generalized least squares problem of finding an estimate $\hat{x}$ given the following model
\begin{equation}
y=Ax+\epsilon,
\end{equation}
where $y$ is a column stacked vector of $n$, $k \times 1$ estimates , each $k\times1$ vector is denoted as $y_i$
\begin{equation}
y=\left[\begin{matrix}y_{1} \\\\\ \vdots \\\\ y_{n}\end{matrix}\right]
\end{equation}
and $A$, our design matrix, is a block matrix, a stack of $n$, $k\times k$ identity matrices
\begin{equation}
A = \left[\begin{matrix} I_{kxk} \\\\ \vdots \\\\ I_{kxk}\end{matrix}\right].
\end{equation}
The size of $A$ is $nk\times k$. Finally, $\epsilon$ is a zero mean Gaussian process with a covariance matrix that overall is unknown but whose block diagonal is provided by the $y_i$ estimated covariance denoted as $E_i$, where $E_i$ is a $k \times k$ matrix. We define the joint covariance matrix $R$ as a block matrix as follows
\begin{equation}
E\left(\epsilon \epsilon^T\right) = R = \left[\begin{matrix} E_{1} & E\left(y_1y_2^T\right) & \hdots & E\left(y_1y_n^T\right) \\\\
E\left(y_2y_1^T\right) & E_2 & \hdots & \vdots \\\\ 
\vdots & \hdots & \ddots  & \vdots \\\\ 
E\left(y_ny_1^T\right) & \hdots & \hdots & E_n
\end{matrix}\right]
\end{equation}
where $E\left(y_iy_j\right)$ is not known when $i \neq j$.

At first this model may look restrictively specific because the design matrix $A$ consists of a column stack of identity matrices. Often in practice measurements must be transformed via a design matrix $A$ which is not a column stack of identity matrices. However, the model above is applicable to data fusion problems where the raw measurements are unavailable. Assume for example that we obtain $n$ estimates of the $k \times 1$ parameters, where each $y_i$ estimate was provided by some black-box processor. This black-box processor might have produced each $y_i$ estimate by fusing some unknown number of raw measurements. Ideally, we might desire to go back to the raw measurements that were used to generate the $n$ estimates. However, this is often not possible; thus, we have the problem presented and discussed in this paper. 

In this paper, we will examine the best linear unbiased estimator (BLUE) for this model under various constructions of the joint unknown covariance matrix $R$. First, we will examine the trivial solution when one assumes that all cross estimate correlations are zero, $E\left(y_iy_j^T\right)_{|i\neq j}= 0$. This approach is sometimes called ellipse convolving. In practice assuming datasets to be completely uncorrelated often produces solution covariance estimates that are too small because the assumption of uncorrelated datasets is often wrong. Some examples of correlation found in datasets thought to be uncorrelated are found in the following sources \cite{pearson1901mathematical, jeffreys1998theory, hampel2011robustsection12,  beck2003systematic, nguyen2011simple}. We will then look at an approach to parameterize the covariance matrix by cross-correlation coefficients. Under this parameterization we will solve for the max-entropy solution as a function of the cross-correlation coefficient for combining two estimates $y_1$ and $y_2$. Finally, we will extend this pairwise approach to $n$ estimates.

\paragraph*{Use of the Terms: Ellipse, Error Integral and Entropy.} The probability density function of the $y_i$ measurement with covariance $E_i$ is 
\begin{equation}
f(y_i)= \left([2 \pi ]^k |E_i|\right)^{-1/2} \exp \left(-\frac{1}{2} (y-y_i)^TE^{-1}(y-y_i) \right).
\end{equation}
We can define the contour where $f(y_i)$ is equal to a constant by finding all whose Mahalanobis distance is the same. In other words, we define the contour as all points $y$ such that
\begin{equation}
\sqrt{(y-y_i)^TE_i^{-1}(y-y_i)} = C,
\end{equation}
where $C$ is a  positive constant.\cite{mahalanobis1936generalized} It is well-known that this contour is the equation of an ellipse for $k=2$ and in general an ellipsoid.  For this reason we use the capital $E$ to denote the covariance and refer to combining measurements as Ellipse Combining. Furthermore, for $k=2$ the area of the ellipse (for $k=3$ the volume of the ellipsoid) defined as the area where
\begin{equation}
\sqrt{(y-y_i)^TE_i^{-1}(y-y_i)} < C
\end{equation}
is proportional to $\sqrt{|E_i|}$.\cite{basu2006data} In this paper we will say that $\sqrt{|E|}$ is proportional to the error integral of the ellipse $E$.
Finally, throughout this paper we will refer to the differential entropy of a continuous random variable simply as the entropy.

\paragraph*{Review of Generalized Least Squares.} For reference and completeness, when the covariance matrix is fully specified and one combines the measurements using a semi-positive definite weight matrix $W$, we have the well-known results \cite{aitken1936iv}
\begin{equation}
\hat{x}=(A^TWA)^{-1}A^TWy,
\end{equation}
\begin{equation} \label{Eq:power}
E\left(\hat{x}\hat{x}^T\right) =  P(W,R) = (A^TWA)^{-1}A^TWRWA(A^TWA)^{-1}.
\end{equation}
The best linear unbiased estimator (BLUE) is found by setting $W=R^{-1}$ and for this specific case Equation \ref{Eq:power} simplifies to

\begin{equation} \label{Eq:powerblue}
E\left(\hat{x}\hat{x}^T\right) =  P(R^{-1},R) = (A^TR^{-1}A)^{-1}.
\end{equation}

As defined so far, the BLUE solution above requires that $R$ be positive definite. We can extend the solution to incorporate $R$ which is only positive semidefinite by setting $W=R^{\dagger}$, where $X^{\dagger}$ is defined as the pseudo-inverse of a positive semidefinite matrix $X$.\cite{gallier2019schur} In the case where $R$ is only positive semidefinite we also require that the solution covariance $P$ be positive definite. For the rest of this paper, for simplicity, we will denote the pseudo-inverse of $R$ as $R^{-1}$ and not $R^{\dagger}$. To summarize, in this paper we require that $R$ be a positive semidefinite such that the resulting solution $P$ is strictly positive definite and $R^{-1}$ will be used to denote the pseudo-inverse of $R$.

\paragraph*{Defining $\alpha$ and $\beta$ to compare models.} We can compare $P(W,R)$ for various $W$ and $R$. 
For example, assume that we have some guess or estimate of the joint covariance denoted $\hat{R}$. We then set  $W^{-1}=\hat{R}$ compute a solution and report a solution covariance via Equation \ref{Eq:power} setting $R=\hat{R}$. Our reported solution covariance will be mismatched if $R \neq \hat{R}$. To describe this mismatch we define
\begin{equation} \label{Eq:alpha}
\alpha(W,R) = \sqrt{\frac{|P(W,R)|}{|P(W,W^{-1})|}}.
\end{equation}
The scalar $\alpha(W,R)$ is the number of times smaller or bigger the error integral of our estimate is compared to what we report.

We also define
\begin{equation} \label{Eq:beta}
\beta(W,R) = \sqrt{\frac{|P(W,R)|}{|P(R^{-1},R)|}}.
\end{equation}
The scalar $\beta(W,R)$ is the number of times bigger the error integral of our estimate is compared to the BLUE solution. The scalar $\beta(W,R) \geq 1$. 

The scalar $\alpha(W,R)$ represents a mismatch between reported covariance and the actual covariance. The scalar $\beta(W,R)$ represents the number of times bigger our solution error integral is compared to what it could be if we fully knew $R$ and applied the BLUE weights (i.e. set $W=R^{-1}$). There are other metrics and methods to compare $P(W,R)$ for various $W$ and $R$; nevertheless, we will use $\alpha$ and $\beta$ as defined above as our metrics for comparison.

\section*{Ellipse Convolving}
We can construct $R$ as $n \times n$ block matrix where each block is a $k \times k$; the size of $R$ is $(nk) \times (nk)$. We denote the $i^{th}$ row and $j^{th}$ column block matrix by $R_{i,j}$. If we assume that all non-diagonal blocks consist of null matrices we have
\begin{equation}
R_{i,i}=E_i,
\end{equation}
\begin{equation}
R_{i,j| j \neq i}=0_{kxk} \:.
\end{equation}
Then we find the BLUE solution is
\begin{equation}
\hat{x}=(A^TR^{-1}A)^{-1}A^TR^{-1}y =  \left(\sum_{i=1}^n E_i^{-1}\right)^{-1} \left(\sum_{i=1}^n E_i^{-1} y_i\right)
\end{equation}
with covariance, 
\begin{equation}
P = (A^TR^{-1}A)^{-1} = \left(\sum_{i=1}^n E_i^{-1}\right)^{-1}.
\end{equation}
For the rest of this paper, we will denote this solution as $\hat{x_c}$ and its covariance as $P_c$ where the subscript $c$ is for ellipse convolving. The underlying assumption is that all provided estimates are independent (and linear). We will consider this the best case scenario and as a result treat $|P_c|$ as a lower bound. In practice, estimates will almost never be completely independent and linear. For example, our estimates might all be off by some unknown bias, which we can treat as a nonzero positive cross correlation between all the estimates. In exploring these other covariance models, we will find that strictly speaking $P_c$ is not always the most optimistic or best case scenario. In other words there are plenty of cases where one can design $R$ to have non-zero positive cross estimate correlation values whose BLUE combine estimate covariance $P$ is such that $|P| < |P_c|$. However, given that cross-estimate correlations are unknown, why select weights that reflect the cross-estimate correlations working together to improve your estimates compared to zero cross correlation? Based on this reasoning, such matrices that theoretically outperform convolving are rejected. Instead, this paper will investigate covariance matrices that result in $|P| \geq |P_c|$.

\section*{Pairwise Ellipse Combining Parameterized by a Scalar Cross-Ellipse Correlation Coefficient}
The previous section showed that ellipse convolving is the best linear unbiased estimator (i.e. BLUE) for combining measurements under the assumption that all cross-ellipse correlation values are zero. This section will investigate BLUE for combining two ellipses whose cross-ellipse correlation is parameterized by a scalar cross-correlation coefficient $r$.  Then we will find the cross correlation coefficient $r_{max}$ which maximizes the entropy of the combined ellipse solution. 

We start by considering the combining of two ellipses with covariance matrices $E_1$ and $E_2$ respectively and a joint system covariance matrix 
\begin{equation}
R = \left[\begin{matrix}E_{1} & r \sqrt{E_{1} E_{2}}\\\\r \left(\sqrt{E_{1} E_{2}}\right)^{T} & E_{2}\end{matrix}\right]
\end{equation} 
where we will define the square root of a matrix $X$ as the principal square root.\cite{higham1987computing} The Schur complement of the $E_2$ block matrix of the matrix $R$ is
\begin{equation} \label{Eq:schur}
R/E_2 = E_1 - r^2  (\sqrt{E_{1} E_{2}}) E_2^{-1}\left(\sqrt{E_{1} E_{2}}\right)^{T}  = (1 - r^2) E_1.
\end{equation}
Clearly the matrix $R/E_2$ is positive semidefinite if and only if $|r| \leq 1$. Or equivalently, $R$ positive semidefinite if and only if $|r| \leq 1$.\cite{gallier2019schur} Also if $|r| = 1$, then $|R|=0$.  Via the block matrix inversion identify \cite{banachiewicz1937berechnung} we find
\begin{equation}
R^{-1} = (1-r^2)^{-1}\left[\begin{matrix}E_{1}^{-1} & -r \sqrt{E_{1}^{-1} E_{2}^{-1}}\\\\-r \left(\sqrt{E_{1}^{-1} E_{2}^{-1}}\right)^{T} & E_{2}^{-1}\end{matrix}\right].
\end{equation}
The entropy of the BLUE solution is 
\begin{equation}
h(P) =  \frac{k}{2} +\frac{k}{2}\ln(2\pi) + \frac{1}{2}\ln \left(|P|\right),
\end{equation}
where $P$ is defined by Equation \ref{Eq:powerblue}.\cite{ahmed1989entropy} Clearly finding the $r_{max}$, the cross-correlation coefficient that maximizes the entropy is equivalent to finding the $r_{max}$ that maximizes $|P|$. For algebraic simplicity, we will examine $|P^{-1}|$,
\begin{equation}
P^{-1} = (1-r^2)^{-1}(S-rZ)
\end{equation}
where we define
\begin{equation}
Z = \sqrt{E_{1}^{-1}E_{2}^{-1}} + \left(\sqrt{E_{1}^{-1}E_{2}^{-1}}\right)^{T},
\end{equation}
\begin{equation}
S = E_{1}^{-1} +E_{2}^{-1}.
\end{equation}
To find the critical points of $|P^{-1}|$ with respect to $r$ we look at
\begin{equation}
\frac{d}{dr} |P^{-1}| = \frac{d}{dr} \left[ (1-r^2)^{-k} |S-rZ| \right]
\end{equation}
and via Jacobi identity\cite{magnus1999matrix} and the derivative chain rule we find
\begin{equation} \label{Eq:dpdr}
\frac{d}{dr} |P^{-1}| = \frac{r^2tr(adj(S-rZ)Z)+2kr|S-rZ|-tr(adj(S-rZ)Z)}{(1-r^2)^{k+1}}
\end{equation}
where $adj(X)$ is the adjugate matrix of $X$ and $tr(X)$ is the trace of $X$.

By looking at the values of $r$ for which $\frac{d}{dr} |P^{-1}|=0$ we can determine which value $r_{max}$ that maximizes $|P|$ for positive coefficient values $0 < r \leq 1$. In doing so, we can find an interval $0 \leq r \leq r_{max} \leq 1$ where the $|P|$ is monotonically increasing with respect to $r$. In other words the entropy of our combine ellipse increases with respect to $r$.
We will restrict ourselves to the range $0 \leq r \leq r_{max}$ when considering alternative algorithms besides convolving to enable ellipse combining. 

It is also useful to compare a mismatch between a weight model based on a predicted $r$ say $r_p$ versus the actual unknown $r$ say we label it $r_n$. For the pairwise case we find that the covariance for the mismatched case given by equation \ref{Eq:power}, now rewritten as a function of $r_p$ and $r_n$ instead of $W$ and $R$ becomes
\begin{equation} \label{Eq:pofrr}
P(r_p,r_n) = (S-r_pZ)^{-1}[(1-2r_pr_n+r_p^2)S+(r_n-2r_p+r_n r_p^2)Z](S-r_pZ)^{-1}.
\end{equation}
We note $P_c = P(0,0)$ and introduce a new constant $P_{max}$ defined as
\begin{equation}
 P_{max} = P(r_{max},r_{max}).
\end{equation}
We can rewrite Equation \ref{Eq:alpha} in terms of $r_p$ and $r_n$ instead of $W$ and $R$, we find
\begin{equation}
\alpha(r_p,r_n)=(1-r_p^2)^{-k/2} \sqrt{\frac{|(1-2r_pr_n+r_p^2)S+(r_n-2r_p+r_n r_p^2)Z| }{|(S-r_pZ)|}}.
\end{equation}
Finally, we can rewrite Equation \ref{Eq:beta} in terms of $r_p$ and $r_n$ as
\begin{equation}
\beta(r_p,r_n)= \sqrt{\frac{|(S-r_nZ)([1-2r_pr_n+r_p^2]S+[r_n-2r_p+r_n r_p^2]Z)|}{(1-r_n^2)^k |(S-r_pZ)|^{2}} }.
\end{equation}
The next two sections provide more specific discussion for 1-D and 2-D location estimates.
\section*{Pairwise Combine: The 1-D case, \bf{\it{k}}=1.} If we have $k=1$ covariance matrices $E_1$ and $E_2$ are scalars; we define
\begin{equation}
E_1 = \sigma_1^2,
\end{equation}
\begin{equation}
E_2 = \sigma_2^2.
\end{equation}
With these definitions, we find
\begin{equation}
Z = \frac{2}{\sigma_1 \sigma_2},
\end{equation}
\begin{equation}
S = \frac{1}{\sigma_1^2} +  \frac{1}{\sigma_2^2}. 
\end{equation}
And equation \ref{Eq:dpdr} simplifies to 
\begin{equation}
\frac{d}{dr} |P^{-1}| = \frac{-Zr^2 +2r(S-rZ)-Z}{(1-r^2)^{2}}.
\end{equation}
Setting the derivative of $|P^{-1}|$ to zero and restricting ourselves to the region $|r| < 1$ we find one critical point and we can easily verify this point is a local maximum of $|P|$, we have
\begin{equation} \label{eq:k1p}
r_{max} = min \left(\frac{\sigma_1}{\sigma_2},\frac{\sigma_2}{\sigma_1} \right).
\end{equation}
We also find
\begin{equation}
\hat{x} = w_1y_1 + w_2y_2
\end{equation}
where 
\begin{equation}
w_1 = \frac{\sigma_2^2-r\sigma_1 \sigma_2}{\sigma_1^2+\sigma_2^2-2r\sigma_1\sigma_2},
\end{equation}
\begin{equation}
w_2 = \frac{\sigma_1^2-r\sigma_1 \sigma_2}{\sigma_1^2+\sigma_2^2-2r\sigma_1\sigma_2}.
\end{equation}
We note that for values $r>r_{max}$ we find either $w_1$ or $w_2$ is less than zero. For the case where the cross-correlation coefficient is unknown and we are computing a weighted average it is counter-intuitive to select a $r$ value that results in a negative weight. Throughout this paper, we avoid this by restricting ourselves to the region $0 \leq r \leq r_{max}$.

\begin{figure}[ht]
\centering
\includegraphics[width=13cm]{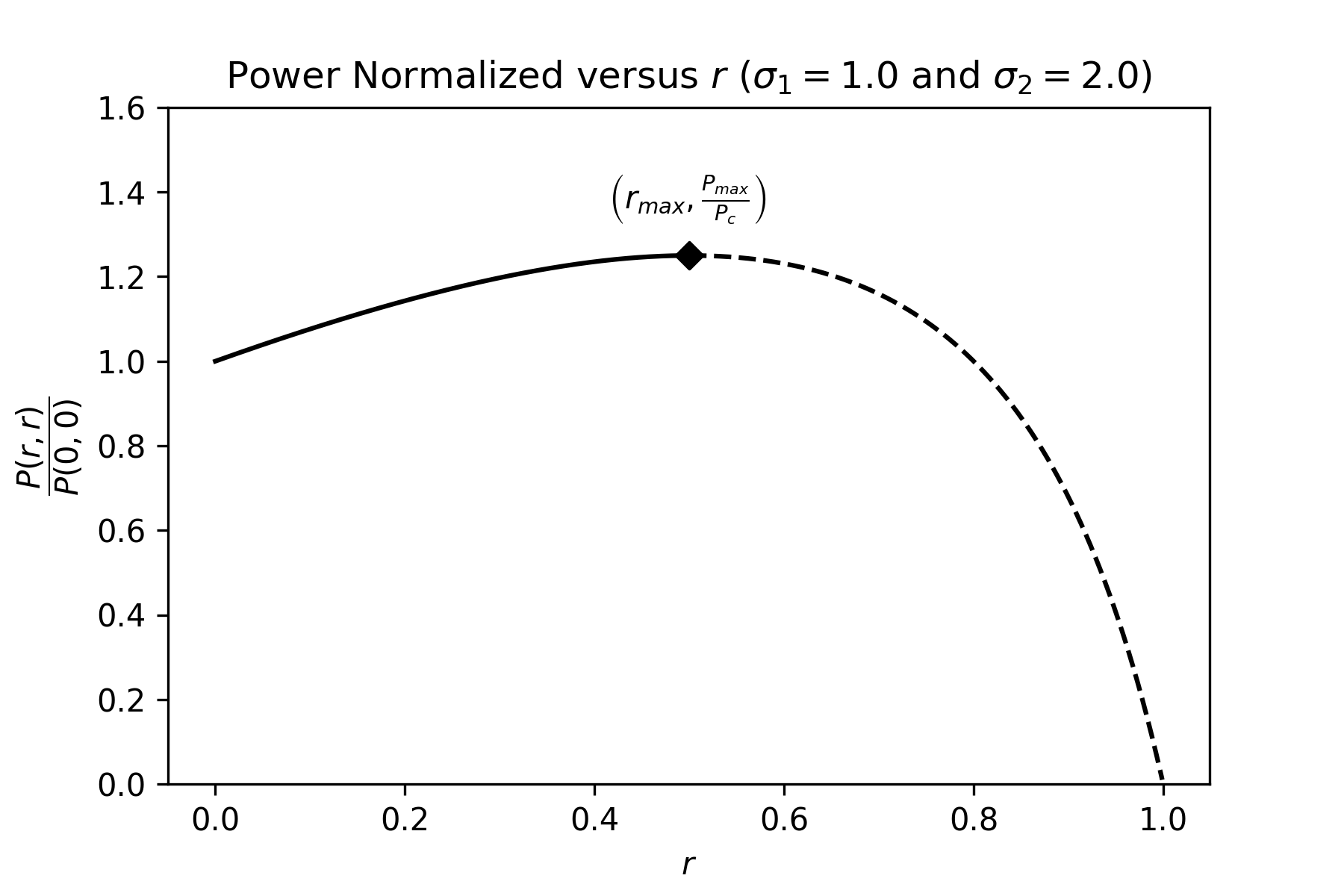}
\caption{Numerical Example of the solution $P$ (defined in equation \ref{Eq:pofrr}) as a function of $r$. The solid line is the interval of $0 \leq r \leq r_{max}$; we recommend only considering $r$ values in this region. The dotted line shows the region where $r_{max} < r \leq 1$.}
\label{fig:pk1}
\end{figure}

\begin{figure}[ht]
\centering
\includegraphics[width=13cm]{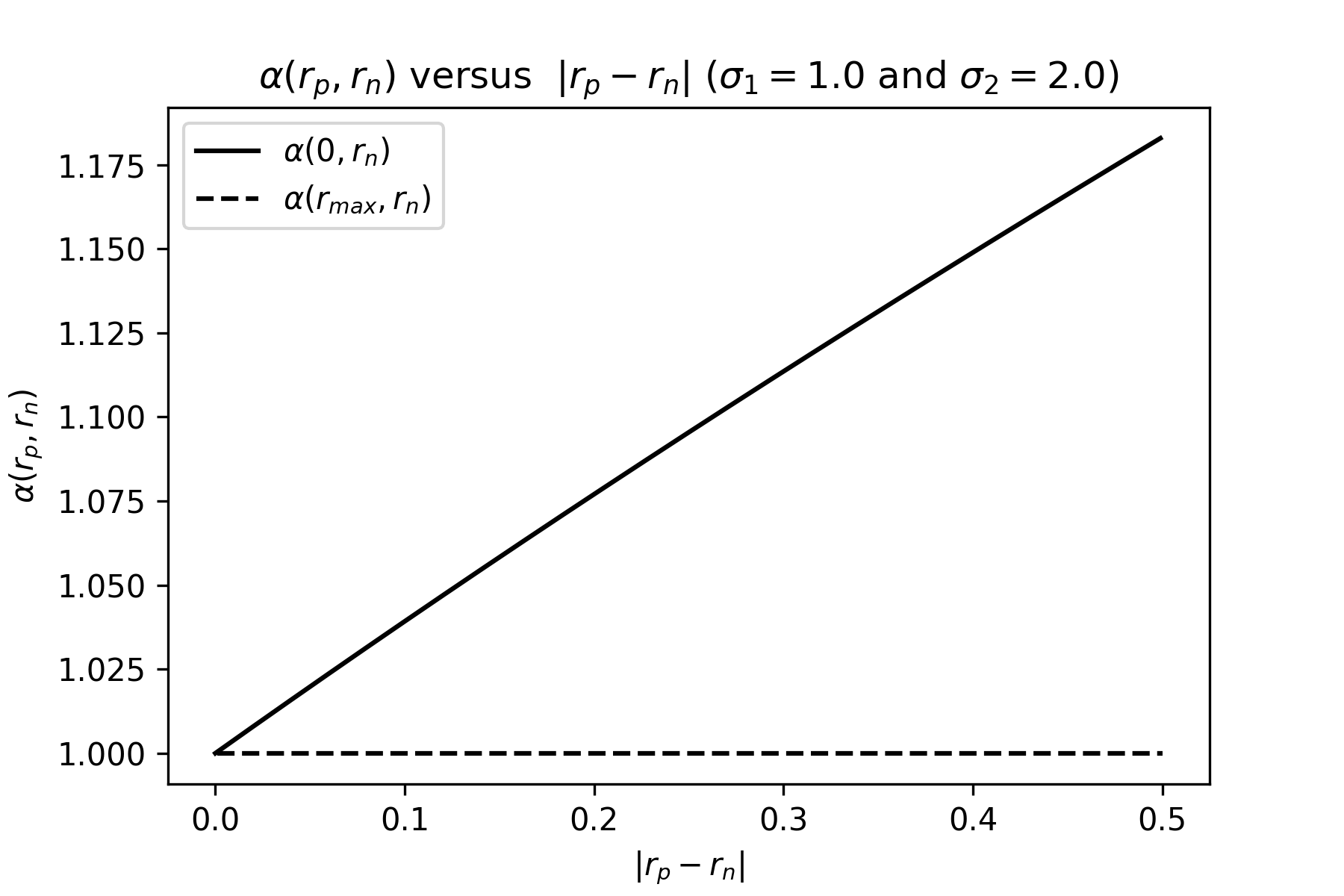}
\caption{Numerical Example of the solution $\alpha(r_p,r_n)$  as a function of $|r_p-r_n|$, where $0 \leq r_n \leq r_{max}$.}
\label{fig:alphak1}
\end{figure}

\begin{figure}[ht]
\centering
\includegraphics[width=13cm]{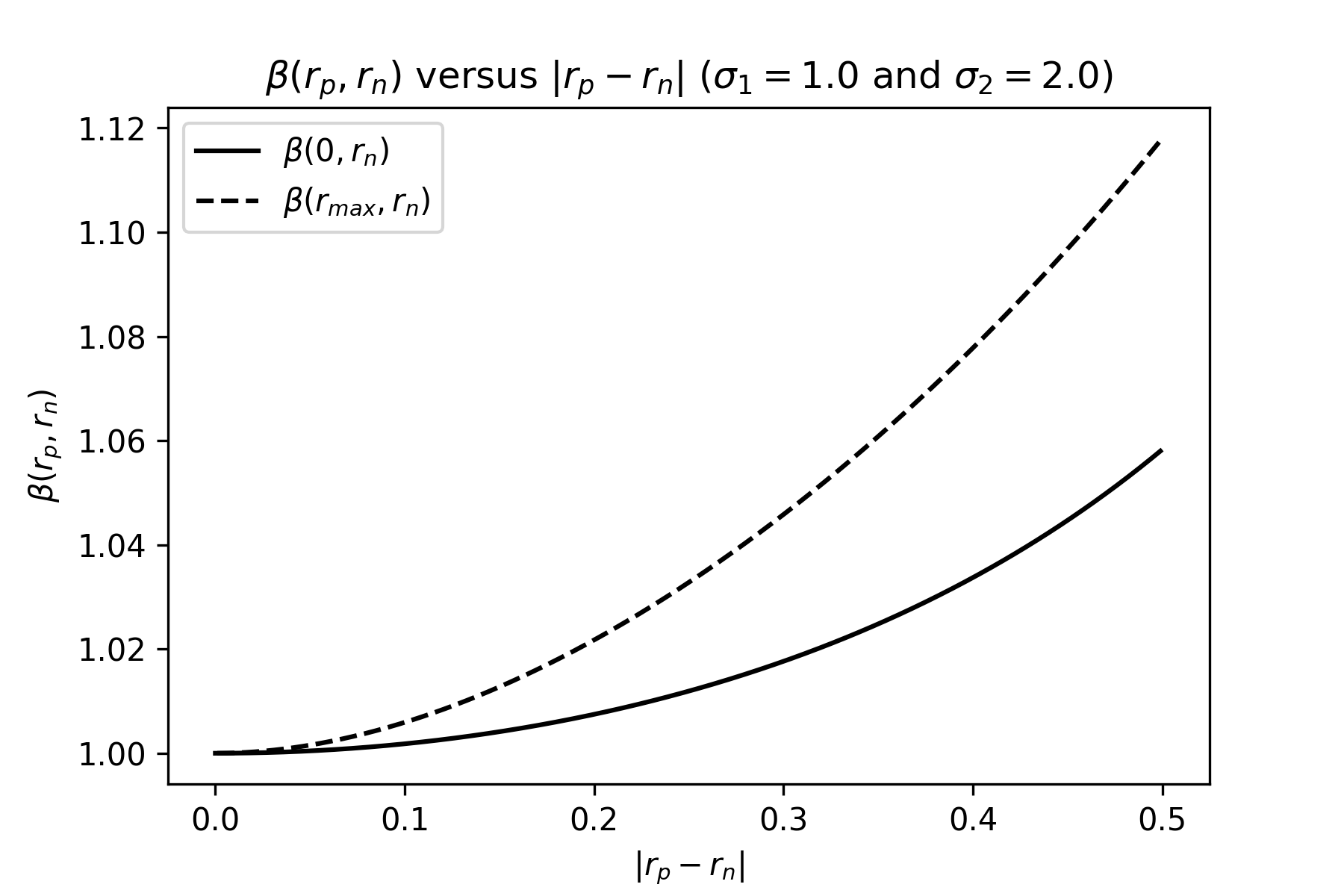}
\caption{Numerical Example of the solution $\beta(r_p,r_n)$  as a function of $|r_p-r_n|$, where $0 \leq r_n \leq r_{max}$.}
\label{fig:betak1}
\end{figure}

For a numerical example, let us consider $\sigma_1 = 1.0$ and $\sigma_2 = 2.0$. Equation \ref{eq:k1p} gives $r_{max} = 0.5$. Figure \ref{fig:pk1} shows the BLUE solution as a function of $r$. Figure \ref{fig:alphak1} shows $\alpha$ for both $r_p=0$ and $r_p=r_{max}$. As $|r_p-r_n|$ increases, $\alpha(0 ,r_n)$ also increases becoming greater than $1$ which means if we convolve and compute the solution covariance assuming the variables are independent and they are not, our solution covariance error integral will be too small. In contrast, as $|r_p-r_n|$ increases, $\alpha(r_{max} ,r_n)$ decreases becoming slightly less than $1$ which means if we compute the combine solution under assume $r=r_{max}$, our computed covariance error integral will be approximately correct and slightly too big if the variables are closer to independent. Finally, Figure \ref{fig:betak1} shows $\beta$ for both $r_p=0$ and $r_p=r_{max}$. From this figure we note that as $|r_p-r_n|$ increases, both $\beta(0 ,r_n)$ and $\beta(r_{max} ,r_n)$ increase. We also note that $\beta(r_{max} ,r_n)$ increases faster than $\beta(0 ,r_n)$ as a function of $|r_p-r_n|$. This means that if we guess $r$ but are off compared to the true $r$, our algorithm will compute a covariance whose error integral is larger than theoretically possible but we will be closer to the lower bound using weights based on $r=0$ (convolving) compared to $r=r_{max}$. 

\section*{Pairwise Combine: The 2-D case, \bf{\it{k}}=2.}
For $2\times 2$ ellipses, we find
\begin{equation}
|P| = \left[\frac{  r^{2} \left|{Z}\right| - \lambda r + \left|{S}\right|}{\left(1 - r^{2}\right)^{2}} \right]^{-1}
\end{equation}
where
\begin{equation}
\lambda = tr \left( S \: adj \left( Z\right)\right).
\end{equation}
We find the $r_{max}$ value will be one of the solutions to the cubic equation,
\begin{equation} \label{Eq:k2cubic}
2|Z|r^3-3\lambda r^2 + (4|S| +2|Z|)r - \lambda = 0.
\end{equation}
We can find solutions to the cubic equation above, throw out solutions where $|r|>1$ and from the remaining solutions find $r_{max}$ and the interval $0 \leq r \leq r_{max}$ where $|P|$ is monotonically increasing with respect to $r$. 

\begin{figure}[ht]
\centering
\includegraphics[width=13cm]{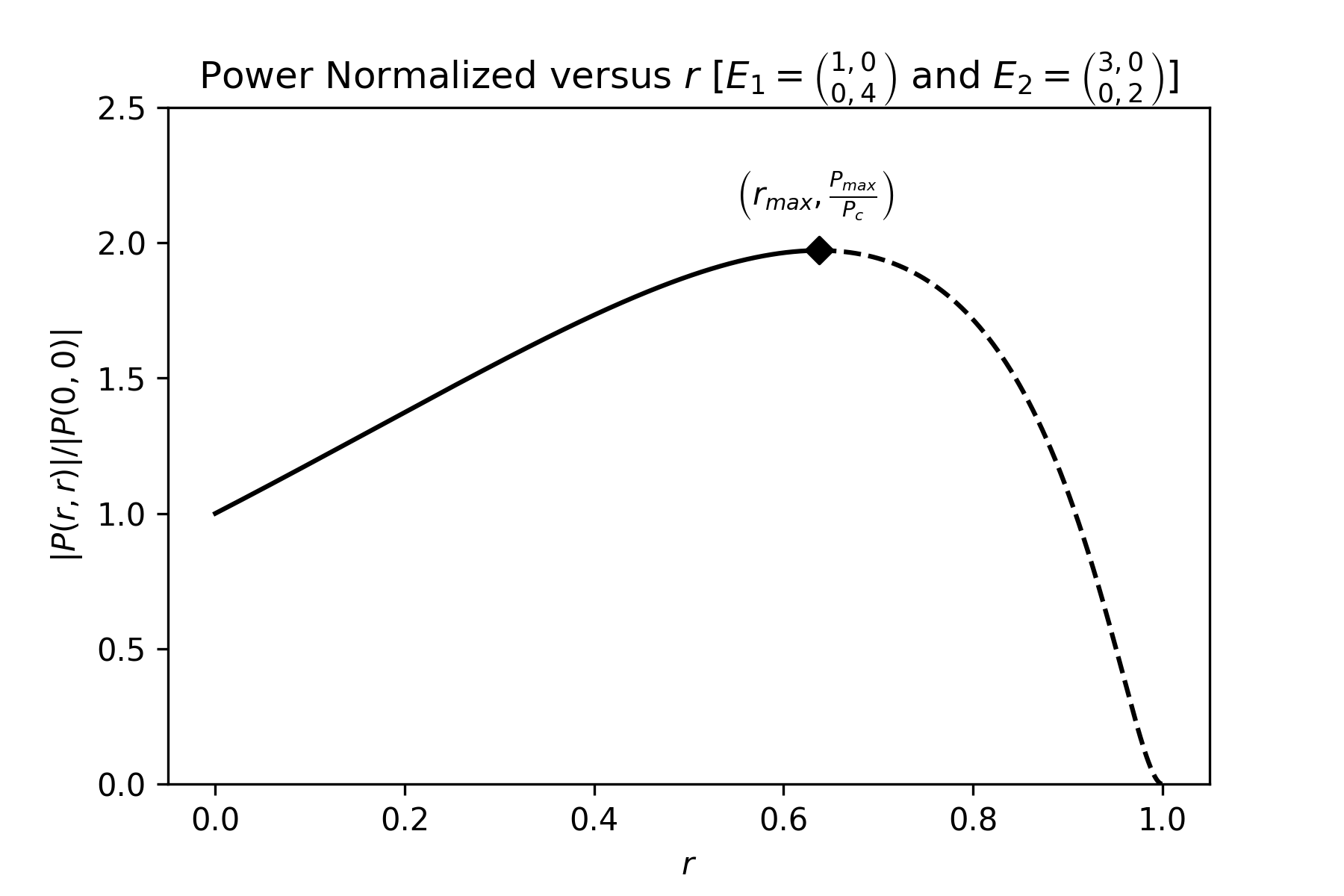}
\caption{Numerical Example of the solution $P$ (defined in equation \ref{Eq:pofrr}) as a function of $r$. The solid line is the interval of $0 \leq r \leq r_{max}$; we recommend only considering $r$ values in this region. The dotted line shows the region where $r_{max} < r \leq 1$.}
\label{fig:pk2}
\end{figure}

\begin{figure}[ht]
\centering
\includegraphics[width=13cm]{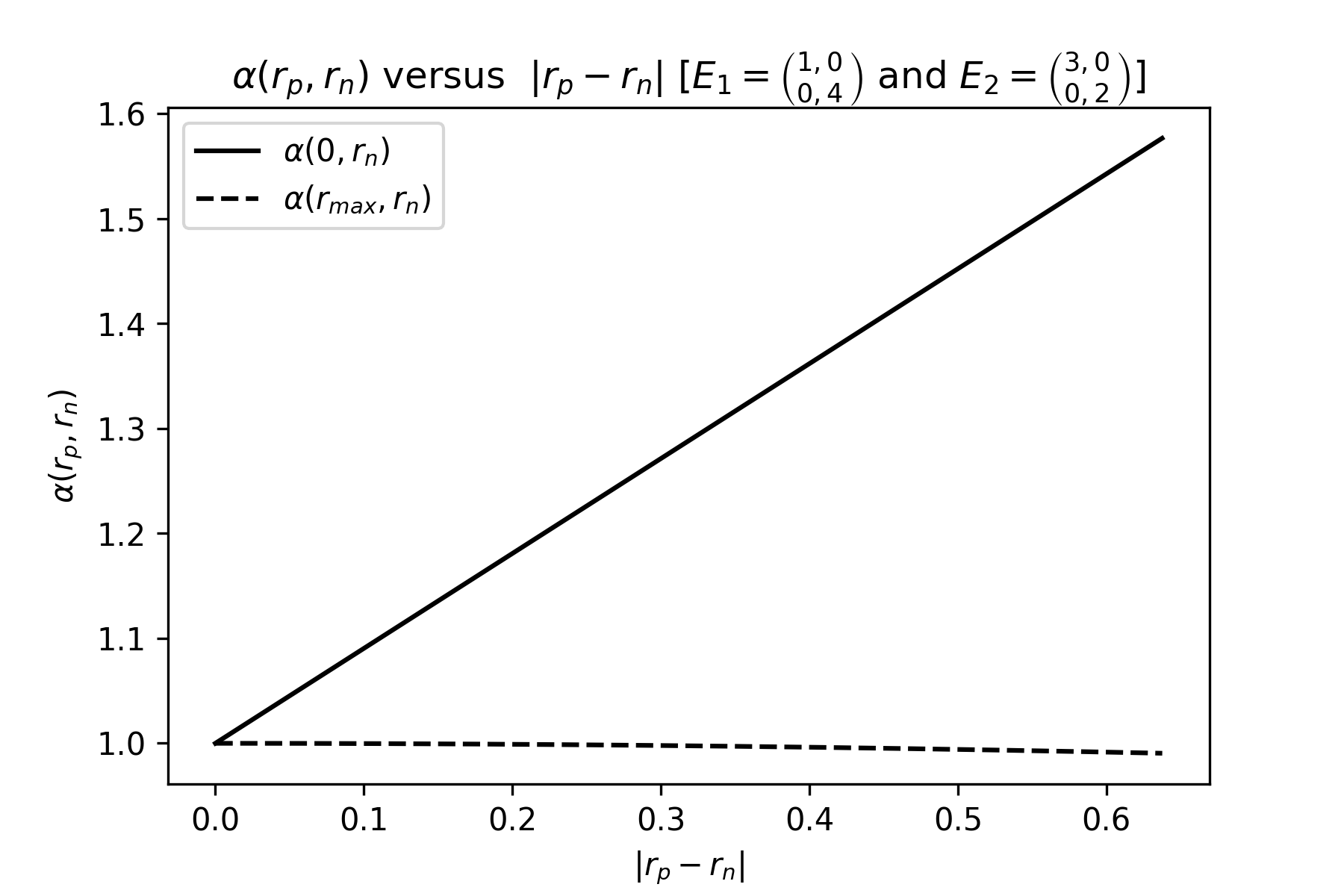}
\caption{Numerical Example of the solution $\alpha(r_p,r_n)$  as a function of $|r_p-r_n|$, where $0 \leq r_n \leq r_{max}$.}
\label{fig:alphak2}
\end{figure}

\begin{figure}[ht]
\centering
\includegraphics[width=13cm]{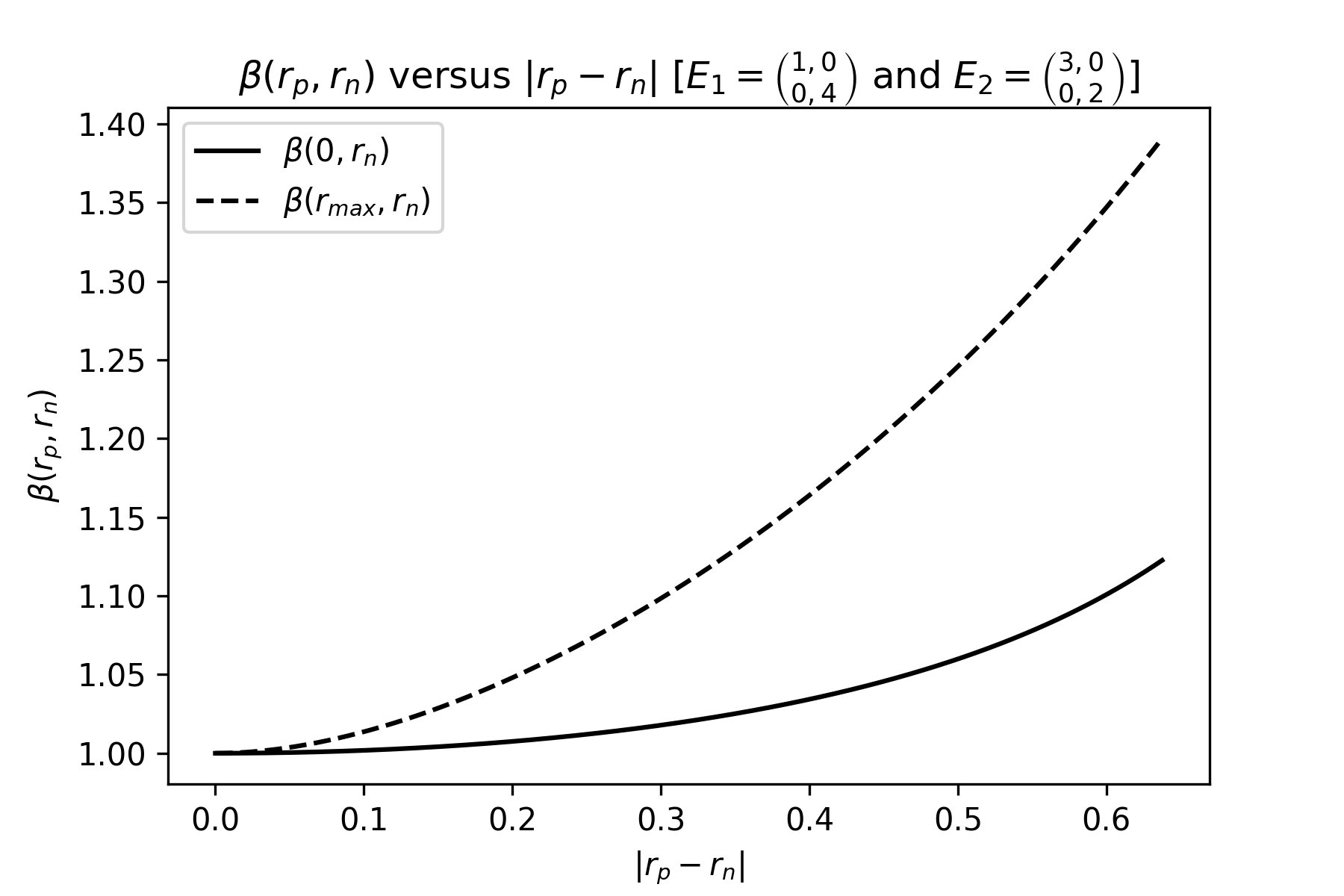}
\caption{Numerical Example of the solution $\beta(r_p,r_n)$  as a function of $|r_p-r_n|$, where $0 \leq r_n \leq r_{max}$.}
\label{fig:betak2}
\end{figure}

For a numerical example let us consider 
\begin{equation}
E_1=\left[\begin{matrix}1 & 0 \\\\0 & 4\end{matrix}\right],
\end{equation}
\begin{equation}
E_2=\left[\begin{matrix}3 & 0 \\\\0 & 2\end{matrix}\right].
\end{equation}
For this example Equation \ref{Eq:k2cubic} becomes
\begin{equation}
\frac{2 \sqrt{6}}{3} r^3 - \frac{3 \sqrt{3}+4\sqrt{2}}{2}r^2 + \frac{2\sqrt{6}+12}{3}r -\frac{2\sqrt{2}}{3}-\frac{\sqrt{3}}{2} = 0
\end{equation}
and we find $r_{max} \approx 0.6376189$. Figure \ref{fig:pk2} shows the BLUE solution as a function of $r$. Figure \ref{fig:alphak2} shows $\alpha$ for both $r_p=0$ and $r_p=r_{max}$. Finally, Figure \ref{fig:betak2} shows $\beta$ for both $r_p=0$ and $r_p=r_{max}$. While the exact values are different, the patterns we observed in the previous section for $k=1$ are the same for $k=2$. Namely, if we construct weights and covariance solution estimation assuming $r=0$ our constructed covariance error integral will be too small if the measurements cross-correlation is non-zero. Conversely, if we construct weights and covariance solution estimation assuming $r=r_{max}$ and $r$ is actually less than $r_{max}$, our constructed covariance error integral will be approximately correct but slightly too large. Additionally, the ratio between the theoretical best error integral and the one computed for a fixed value of $|r_p-r_n|$ is smaller when $r_p=0$ compared to $r_p=r_{max}$.

\section*{Extending Ellipse Combining to $n$ Ellipses}
In general, we will have more than two ellipses we wish to combine. The extension is straightforward. We parameterize the joint system covariance matrix by correlation coefficient between every pair of ellipses. We consider $n$ ellipses with covariance matrices $E_1,E_2,...,E_j,...,E_n$ all of size $k\times k$. Then, we define the joint covariance matrix parameterized by pairwise correlation coefficients as a block matrix consisting of $n\times n$ blocks. Each block is a square matrix equal in size to the individual ellipse covariance. The $i$ row, $j$ column sub-block is defined as 
\begin{equation}
R_{i,j|i<j} = r_{i,j} \sqrt{E_iE_j},
\end{equation}
\begin{equation}
R_{i,i} = E_i,
\end{equation}
\begin{equation}
R_{i,j|i>j} =  r_{i,j}\left(\sqrt{E_{i} E_{j}}\right)^{T}.
\end{equation}
Or in other words we define
\begin{equation}
R(\overrightarrow{r}) = \left[\begin{matrix} E_{1} & r_{1,2} \sqrt{E_1E_2} & \hdots & r_{1,n} \sqrt{E_1E_n} \\\\
r_{1,2}\left(\sqrt{E_{1} E_{2}}\right)^{T} & E_2 & \hdots & \vdots \\\\ 
\vdots & \hdots & \ddots  & \vdots \\\\ 
r_{1,n}\left(\sqrt{E_{1} E_{n}}\right)^{T} & \hdots & \hdots & E_n
\end{matrix}\right]
\end{equation}
where $\overrightarrow{r}$ is a vector of $\frac{n^2-n}{2}$ correlation coefficients. The values of $\overrightarrow{r}$ are constrained to values such that $R(\overrightarrow{r})$ is positive semidefinite.

Over all possible valid $\overrightarrow{r}$ we define $\overrightarrow{r}_{max}$ as the vector that maximizes $P(R^{-1}(\overrightarrow{r}),R(\overrightarrow{r}))$ as defined in Equation \ref{Eq:power}. While identifying and developing algorithms to efficiently find $\overrightarrow{r}_{max}$ is left for future work, we will define an initial guess of $\overrightarrow{r}_{max}$ based on the pairwise solution. We define $\overrightarrow{r}_{pm}$ as the vector where each $r_{i,j}$ is set by considering the pairwise problem involving just $E_i$ and $E_j$ and finding the coefficient that maximizes this pairwise problem as defined in the previous sections. The subscript $pm$ in $\overrightarrow{r}_{pm}$ denotes the pairwise max. We can use $\overrightarrow{r}_{pm}$ as an initial guess of $\overrightarrow{r}_{max}$ or even as a computationally inexpensive substitute. We conjecture that the matrix $R(\overrightarrow{r}_{pm})$ is always positive semidefinite, or in other words that $\overrightarrow{r}_{pm}$ is a valid choice of $\overrightarrow{r}$. 

We define $R_{max}$ and $R_{pm}$ as
\begin{equation} \label{Eq:rmax}
R_{max} = R(\overrightarrow{r}_{max}),
\end{equation}
and
\begin{equation} \label{Eq:rpm}
R_{pm} = R(\overrightarrow{r}_{pm}).
\end{equation}
$R_{max}$ is the covariance matrix that results in the blue solution whose uncertainty has the largest entropy, and $R_{pm}$ is an initial guess or practical estimate of $R_{max}$. It is noted that in the specific case of $k=1$, $R_{pm}=R_{max}$.

\section*{Three Alternatives to Convolving}
The convolving algorithm and resulting combined covariance are based on an assumption of independent and zero-mean measurements. If we have reason to believe that our observations are correlated or biased we can adopt a different strategy. While there are numerous possible strategies, in this section we will discuss three.
\paragraph*{Max Entropy BLUE.} Given $n$ observations ${y_i}$ and their respective covariance matrices ${E_i}$, we build the max entropy parameterized covariance matrix $R_{max}$ as defined by Equation \ref{Eq:rmax}. Then, compute the blue solution. The resulting combined estimate and reported covariance are given by
\begin{equation}
\hat{x}=(A^TR_{max}^{-1}A)^{-1}A^TR_{max}^{-1}y,
\end{equation}
\begin{equation}
P_{max} = (A^{T}R_{max}^{-1}A)^{-1}.
\end{equation}
This algorithm assumes the cross-observation correlations are such to make our estimation as uncertain as possible. Compared to the next proposed algorithm, this algorithm is computationally expensive especially as the number of measurements $n$ increases. Computing $\hat{x}$ directly involves both finding $R_{max}$ and taking the inverse of $R_{max}$, a matrix of size $(nk) \times (nk)$. We note, one could also use $R_{pm}$ in place of $R_{max}$ in order to avoid the computational cost of finding $R_{max}$.
\paragraph*{Assume Max Entropy to Report Covariance; Convolve to estimate solution.} Given $n$ observations ${y_i}$ and their respective covariance matrices ${E_i}$, we build the max entropy parameterized covariance matrix $R_{max}$ as defined by Equation \ref{Eq:rmax}. Next, we estimate the solution assuming the measurements are uncorrelated (i.e. convolve) and report the combined estimate covariance under the max entropy assumption. We define
\begin{equation}
P_c = \left(\sum_{i=1}^n E_i^{-1}\right)^{-1},
\end{equation}
and 
\begin{equation}
P_r = P_c \left(\sum_{i=1}^n\sum_{j=i+1}^n r_{i,j,max}  \left[\sqrt{E_{i}^{-1}E_{j}^{-1}} + \sqrt{E_{j}^{-1}E_{i}^{-1}}\right]\right) P_c.
\end{equation}
Subsequently, we find that in this algorithm $\hat{x}$ and reported covariance $P$ as follows, 
\begin{equation}
\hat{x}=  P_c \left(\sum_{i=1}^n E_i^{-1} y_i\right)
\end{equation}
\begin{equation}
P= P_c + P_r.
\end{equation}
This algorithm is convolving with the only difference of a modified reported covariance. The modified covariance matrix is a conservative estimate which assumes the max entropy joint covariance $R_{max}$. In contrast to the previous algorithm $R_{max}$ is never inverted. Additionally, if one assumes that new measurements are coming in over time, it is easy to write this algorithm as an update to previous solutions with each new estimate. 
\paragraph*{Max Entropy provided a sum of covariances per measurement.} This approach requires that each observation covariance is specified as a linear sum of covariance to include an independent portion and some number, say $m$ correlated or bias portions,
\begin{equation}
E_i = B_0 + \left(\sum_{a=1}^m B_a \right). 
\end{equation}
Then we define the joint covariance model as
\begin{equation}
R = R_0(\overrightarrow{0}) + \left(\sum_{a=1}^m R_{a}(\overrightarrow{r_a}) \right),
\end{equation}
where  $\overrightarrow{r_a}$ is a vector of correlation coefficients for the $a^{th}$ covariance matrix in the sum. We define the number of covariance terms $m$ and their respective $\overrightarrow{r_a}$ based on the specific model we wish to adopt. We have defined $B_0$ to be the portion of each  covariance that is completely independent and therefore the  joint covariance $R_0$ is just a block diagonal matrix. After constructing $R$, the algorithm performs the best linear unbiased estimator. For example, we might have a bias covariance due to a common instrument used to make the observations so for all measurements that come from that instrument we select the max entropy coefficient within $\overrightarrow{r_a}$. While if observations come from another instrument, we set the appropriate $\overrightarrow{r_a}$ term to zero.
This approach requires that the measurements provided come with more than just a single covariance matrix to enable one to construct $R$. One could make these models very complex with a large $m$ and different rules for building each  $\overrightarrow{r_a}$, but doing so would require that each measurement be stored with all of these covariance terms. Specific approaches will depend on each application; we present one simple version of this approach which leverages the time of a measurement in the construction of $R$.

A simple version of this approach which enables weight models that emphasize time diversity is a case where  we set $m=1$, or assume that each observation is broken into a independent covariance term, $B_0$, and a bias or correlated bias term, $B_1$. Additionally, assume each measurement comes with a measurement time, $t_i$. Then for the bias vector $\overrightarrow{r_1}$ compute the pairwise component as
\begin{equation}
r_{i,j} = r_{i,j,max} \exp(-\gamma|t_i-t_j|),
\end{equation}
where $\gamma$ is a positive real constant. This model reflects a  system where measurements close in time are correlated, but this correlation decreases as the time between measurements increases.

\section*{Conclusion}
Measurement combining via convolving assumes provided measurements are uncorrelated. When the measurements are not uncorrelated the reported covariance that comes from convolving is often too small. We show an approach to construct joint measurement covariance matrices parameterized by measurement pair-wise correlation coefficients. Finally based on this approach, we present three alternatives to convolving when the full-joint covariance is not provided and must be constructed. 

In practice, linear uncorrelated location estimates almost never exist. When we assume estimates are uncorrelated, the combined estimated covariance determinant shrinks to zero as the number of measurements increase. In contrast, this is not true for the alternative algorithms presented in this paper. In real world applications, believing we can achieve arbitrarily precise estimates by increasing the number of measurements is naive (for further discussion and examples see \cite{hampel2011robustsection81, youden1972enduring, mosteller1977data}).  We therefore strongly recommend algorithms like the ones discussed in this paper in place of convolving. While this paper focused on cases where the mapping of a single measurement vector to the solution space is the identity matrix, the concepts and ideas are extendable to other design matrices. 

Lastly, we note that all algorithms presented are linear. It is well-known that linear estimates are particularly sensitive to outliers or data contamination. Appropriate modifications to the presented algorithms to make them robust are recommended.

\bibliographystyle{plainnat}
\bibliography{EllipseCombiningWithUnknownCrossEllipseCorrelations}

\begin{thebibliography}{16}
\providecommand{\natexlab}[1]{#1}
\providecommand{\url}[1]{\texttt{#1}}
\expandafter\ifx\csname urlstyle\endcsname\relax
  \providecommand{\doi}[1]{doi: #1}\else
  \providecommand{\doi}{doi: \begingroup \urlstyle{rm}\Url}\fi

\bibitem[Ahmed and Gokhale(1989)]{ahmed1989entropy}
Nabil~Ali Ahmed and DV~Gokhale.
\newblock Entropy expressions and their estimators for multivariate
  distributions.
\newblock \emph{IEEE Transactions on Information Theory}, 35\penalty0
  (3):\penalty0 688--692, 1989.

\bibitem[Aitken(1936)]{aitken1936iv}
Alexander~C Aitken.
\newblock Iv.—on least squares and linear combination of observations.
\newblock \emph{Proceedings of the Royal Society of Edinburgh}, 55:\penalty0
  42--48, 1936.

\bibitem[Banachiewicz(1937)]{banachiewicz1937berechnung}
Tadeusz Banachiewicz.
\newblock Zur berechnung der determinanten, wie auch der inversen und zur
  darauf basierten auflosung der systeme linearer gleichungen.
\newblock \emph{Acta Astronom. Ser. C}, 3:\penalty0 41--67, 1937.

\bibitem[Basu and Ho(2006)]{basu2006data}
Mitra Basu and Tin~Kam Ho.
\newblock \emph{Data complexity in pattern recognition}, page 181.
\newblock Springer Science \& Business Media, 2006.

\bibitem[Beck et~al.(2003)Beck, Shukurov, Sokoloff, and
  Wielebinski]{beck2003systematic}
Rainer Beck, Anvar Shukurov, Dmitry Sokoloff, and Richard Wielebinski.
\newblock Systematic bias in interstellar magnetic field estimates.
\newblock \emph{Astronomy \& Astrophysics}, 411\penalty0 (2):\penalty0 99--107,
  2003.

\bibitem[Gallier(2019)]{gallier2019schur}
Jean Gallier.
\newblock The schur complement and symmetric positive semidefinite (and
  definite) matrices.
\newblock 2019.

\bibitem[Hampel et~al.(2011{\natexlab{a}})Hampel, Ronchetti, Rousseeuw, and
  Stahel]{hampel2011robustsection12}
Frank~R Hampel, Elvezio~M Ronchetti, Peter~J Rousseeuw, and Werner~A Stahel.
\newblock \emph{Robust statistics: the approach based on influence functions},
  volume 196, pages 22--25.
\newblock John Wiley \& Sons, 2011{\natexlab{a}}.

\bibitem[Hampel et~al.(2011{\natexlab{b}})Hampel, Ronchetti, Rousseeuw, and
  Stahel]{hampel2011robustsection81}
Frank~R Hampel, Elvezio~M Ronchetti, Peter~J Rousseeuw, and Werner~A Stahel.
\newblock \emph{Robust statistics: the approach based on influence functions},
  volume 196, pages 387--397.
\newblock John Wiley \& Sons, 2011{\natexlab{b}}.

\bibitem[Higham(1987)]{higham1987computing}
Nicholas~J Higham.
\newblock Computing real square roots of a real matrix.
\newblock \emph{Linear Algebra and its applications}, 88:\penalty0 405--430,
  1987.

\bibitem[Jeffreys(1998)]{jeffreys1998theory}
Harold Jeffreys.
\newblock \emph{The theory of probability}.
\newblock OUP Oxford, 1998.

\bibitem[Magnus and Neudecker(1999)]{magnus1999matrix}
J.R. Magnus and H.~Neudecker.
\newblock \emph{Matrix Differential Calculus with Applications in Statistics
  and Econometrics}, pages 149--150.
\newblock Wiley, 1999.
\newblock ISBN 9780471986331.

\bibitem[Mahalanobis(1936)]{mahalanobis1936generalized}
Prasanta~Chandra Mahalanobis.
\newblock On the generalized distance in statistics.
\newblock National Institute of Science of India, 1936.

\bibitem[Mosteller et~al.(1977)Mosteller, Tukey, et~al.]{mosteller1977data}
Frederick Mosteller, John~Wilder Tukey, et~al.
\newblock \emph{Data analysis and regression: a second course in statistics},
  chapter Hunting Out the Real Uncertainty.
\newblock 1977.

\bibitem[Nguyen and Jiang(2011)]{nguyen2011simple}
Thuan Nguyen and Jiming Jiang.
\newblock Simple estimation of hidden correlation in repeated measures.
\newblock \emph{Statistics in medicine}, 30\penalty0 (29):\penalty0 3403--3415,
  2011.

\bibitem[Pearson(1901)]{pearson1901mathematical}
Karl Pearson.
\newblock On the mathematical theory of errors of judgement with special
  reference to the personal equation.
\newblock \emph{Proceedings of the Royal Society of London}, 68\penalty0
  (442-450):\penalty0 369--372, 1901.

\bibitem[Youden(1972)]{youden1972enduring}
WJ~Youden.
\newblock Enduring values.
\newblock \emph{Technometrics}, 14\penalty0 (1):\penalty0 1--11, 1972.

\end{thebibliography}

\end{document}